%Paper: dg-ga/9501006
%From: Alexander Reznikov <simplex@MATH.HUJI.AC.IL>
%Date: Mon, 30 Jan 1995 11:23:36 +0200

\input vanilla.sty
\magnification 1200
\baselineskip 18pt
\input definiti.tex
\input mathchar.tex

\define\pmf{\par\medpagebreak\flushpar}

\pmf
\title Simpson's Theory and Superrigidity of Complex Hyperbolic Lattices
\endtitle
\author
Alexander Reznikov
\endauthor

{\bf Abstract} \ We attack a conjecture of J. Rogawski:  any cocompact lattice
in $S U (2, 1)$ for which the ball quotient $X = B^2 / \Gamma$ satisfies $b_1
(X) = 0$ and
$H^{1, 1} (X) \cap H^2 (X, \bbq) \approx \bbq$ is arithmetic.  We prove the
Archimedian suprerigidity for representation of $\Gamma$ is $S L (3, \bbc)$.

\title
Th\'eorie de Simpson et superrigidit\'e des r\'eseaux hyperboliques complexes
\endtitle

{\bf R\'esum\'e} Soit $\Gamma \subset S U (2,1)$ un reseau cocompact et soit $X
= B^2 / \Gamma$.  Nous preuvons: si $b_1 (X) = 0$ et $H^{1,1} (X) \cap
H^2 (X, \bbq) \approx \bbq$ allors tous les representations $\rho$ de $\Gamma$
dans $S L (3, \bbc)$ sont conjugu\'e \`a le repr\'esentation naturelle ou la
fermeture de Zariski de l'image $p (\Gamma)$ est compacte.

{\bf Version fran\c caise abr\'eg\'e\'e} - Le th\'eor\`eme classique
de Margulis dit que tous les
r\'eseaux dans les groupes de Lie semi-simples sont superrigides.  Ceci a et\'e
generalis\'e par Corlette [C] \`a la superrigidite
des r\'eseaux quaternioniques et de Cayley.  D'autre part, Johnson et
Millson ont montr\'e q\'u il existait des deformations
des r\'eseaux cocompact dans $SO (n, 1)$ si on regarde $SO (n,1)$
comme plong\'e dans $SO (n +1, 1)$.

C'est une question d'un int\'er\^et fondamental de savoir si les
r\'eseaux hyperboliques complexes sont superrigides.

Dans cet article, nons considerons la question suivante de J. Rogawski.

\underbar{Hypoth\'ese}  Soit $X = B^2 / \Gamma, \Gamma \subset S U (2, 1)$
une surface hyperbolique complexe compacte.
Supposons $b_1 (X) = 0$ et $H^{1,1} (X) \cap H^2 (X, \bbq) = \bbq$.  Allors
$\Gamma$
est arithmetique et provient d'une alg\'ebre avec division $E | \bbq$
de rang 3 avec une involution.

Observons que pour tous les r\'eseux provenant d'alg\'ebres avec division, on a
effectivement $b_1 (X) = 0$
 et $H^{1,1} (X) \cap H^2 (X, \bbq) = \bbq$ [Rog].

Soit $\ell$ un fibr\'e lin\'eaire tautologiue sur $X$  [Re].  La condition
$H^{1,1} (X) \cap H^2 (X, \bbq) = \bbq$ dit que $[\ell] = k \cdot$
g\'en\'erateur dans $Pic (X) / tors \approx \bbz$.

Le r\'esultat principlal de cet article prouve la superrigidit\'e des
representations de $\Gamma$ dans $S L (3, \bbc)$ dans le cas $k =1$.

\underbar{Theorem\'e principal}  Soit $X = B^2/ \Gamma$ et supposons que $b_1
(X) = 0$ et
$H^{1,1} (X) \cap H^2 (X, \bbq) = \bbq$.  Soit $[\ell]$ un g\'en\'erateur de
$Pic (X) / tors \approx \bbz$.
Si $\rho$ est une repr\' sentation de $\Gamma$ dans $S L (3, \bbc)$ alors soit
$\rho$ est conjugu\'e \`a la repr\'esentation naturelle de $\Gamma$, soit la
fermeture de
Zariski de l'image $\rho (\Gamma)$ est compacte.

Je vaudrais remercier Ron Livn\'e,Jon Rogawski, et Carlos Simpson
 pour beacoup de discussiones int\'eresantes.  Je
vaudrais aussi remercier Marina Ville et Lucy Katz pour son aide essentielle
\`a la pr\'eparation de cet article.
\par
\newpage
\par
\title
Simpson's Theory and Superrigidity of Complex Hyperbolic Lattices
\endtitle
\author
Alexander Reznikov
\endauthor

{\bf 0} \underbar{Main Theorem}  The classical theorem of Margulis establishes
the
superrigidity of lattices in semisimple Lie groups of rank $\ge 2$.  The
work of Corlette [C] extended this to (Archimedian) superrigidity of uniform
quaternionic and Cayley lattices.  On the other hand, by Johnson and
Millson [JM] some uniform lattices in
$SO (n, 1)$ admit deformations as mapped to $SO (n +1, 1)$.

It is therefore of fundamental interest to study to what extent the
complex hyperbolic lattices are superrigid.  Since there are nontrivial
holomorphic
maps between different ball quotients [DM] one should confine oneself's look to
lattices (or manifolds) ``minimal'' in some sense.

The present note addresses the following conjecture of Jon Rogawski.

\underbar{Conjecture}.  Let $X = B^2 / \Gamma, \Gamma \subset S U (2, 1)$ be a
compact ball quotient.
Suppose $b_1 (X) = 0$ and $H^{1,1} (X) \cap H^2 (X, \bbq) \approx \bbq$.  Then
$\Gamma$ is arithmetic and comes from a division algebra $E | \bbq$ of rank 3
with an involution.

Observe that for all lattices coming from division algebras, indeed $b_1 (X) =
0$ and $H^{1,1} (X) \cap H^2 (X, \bbq) \approx
\bbq$ [Rog].

Let $\ell$ be the tautological line bundle over $X$ [Re].  Since
$Pic (X) / tors \approx \bbz$,
we have $[ \ell] = k \cdot$ generator for some $k \in \bbz$.

The main result of the paper establishes the superrigidity of representations
of $\Gamma$ in $SL (3, \bbc)$ for $\Gamma$
yielding $k =1$ as follows.

\underbar{Main Theorem}. Let $X = B^2 / \Gamma$ and suppose $b_1 (X) = 0$
and $H^{1,1} (X) \cap H^2 (X, \bbq) \approx Q$.  If
[$\ell$] generates $Pic (X) / tors \approx \bbz$, then any representation
of $\Gamma = \pi_1 (X)$ in $S L (3, \bbc)$ is either conjugate to the
natural representation up to the twist by a character, or has a compact
Zariski closure.

One hopes, that, applying methods of [GS]  one is able to prove
the $p$-adic superrigidity and to settle Rogawski's conjecture.

I wish to thank Ron Livne, Jon Rogawski and Carlos Simpson for stimulating
discussions.

{\bf 1.} \underbar{Computations of Higgs bundles}.  We admit a knowledge
of Simpson's theory [S1].  Let $X$ be as above and let $\rho_0 :
\Gamma \to PS U (n, 1)$ be the natural
representation. Then the corresponding Higgs bundle is as follows [Re].  Take
$E = T X \otimes \ell \oplus \ell$ as a holomorphis bundle and define
$\theta \in H^0 (\Omega^1 \otimes End (E))$ by $\pmatrix
0 &0 \\
0 &1 \endpmatrix$.

In view of the Simpson's theory, for proving the Main Theorem one needs to show
that any complex variation of Hodge structure [S1] of type (2,1) over $X$
is as above.  Indeed, any representation is deformable to one, corresponding to
a variation
of Hodge structure [S2], and the natural representation is rigid [W].

So let $F = ( \xi \oplus \eta, \theta)$ be a variation of complex Hodge
structure, rank $\xi = 2$, rank $\eta =1, \theta \in H^0
(\Omega^1 (X) \oplus
\text{Hom} (\eta, \xi)) \approx H^0 (\text{Hom} (T X \otimes \xi, \eta))$.

\proclaim{1.2. \underbar{Lemma}}  Let $\lambda, \mu$ be rank two bundles over
$X$
and let $f \in H^0 (\text{Hom} (\lambda, \mu)), f \neq 0$.  Then either
rank $f \le 1$ everywhere or
$$ (c_1 (\mu), [\omega] ) \ge (c_1 (\lambda), [ \omega]) $$
with the equality iff $\lambda \approx \mu$, and rank $f = 2$ everywhere.
Here $[\omega]$ is the K\"ahler class.
\endproclaim

\demo{Proof}  Consider $\wedge^2 f: \wedge^2 \lambda \to \wedge^2 \mu$.
If $\wedge^2 f \neq 0$,
then $\wedge^2 \mu \oplus (\wedge^2 \lambda)^{-1}$ has a nontrivial
holomorphic section, whose zero locus is an effective divisor, so
$(c_1 (\wedge^2 \mu \otimes (\wedge^2 (\lambda))^{-1} ), [ \omega] ) \ge 0$ and
the
equality implies $\wedge^2 f$ is an isomorphism.

\demo{2. \underbar{Proof of the Main Theorem}}

Let $F = (\xi \oplus \eta, \theta)$ be as above.

\underbar{Case 1} Rank $\theta = 2$ somewhere.

Applying the lemma, we get
$$ (c_1 (T X \otimes \xi), [ \omega] ) \le (c_1 (\eta), [ \omega] ) $$
Now, $[\omega] \sim [\ell]$ since $X$ is hyperbolic, and $c_1 (TX) = - 3 [
\ell]$ in $H^2 (X, \bbr)$, so
$$ (c_1 (\eta) - 2 c_1 (\xi), [ \ell] ) \le 3 [ \ell ]^2. $$
On the other hand, since $\xi \oplus \eta$ is a deformation of the flat
bundle, $c_1 (\xi \oplus \eta) = 0$, i.e. $c_1 (\xi) = - c_1 (\eta)$, so
$$ (c_1 (\xi), [ \ell] ) \le [ \ell]^2. \tag * $$
Since $F$ is $\theta$-stable [S1], $(c_1 (\eta), [ \ell] ) < 0$, so
$(c_1 (\eta), [ \ell]) > 0$.  This leaves the only possibility
$(c_1 (\xi), [ \ell]) = [ \ell]^2$,
because $[ \ell]$ generates $Pic (X) / tors$.  So $\xi = \ell \otimes \alpha$,
where $\alpha$ is a linear unitary flat bundle, corresponding to
$Pic (X)/tors \approx H^{tors}_1 (X, \bbz)$.
(recall that $b_1 (X) = 0$).  Moreover, since (*) becomes
an equality, we get by lemma above $\eta \approx T X \otimes \xi$ and
$\theta$
takes the form $\pmatrix 0 &0 \\
0 &1 \endpmatrix$.  Hence $F \approx (T X \otimes \ell \oplus \ell) \otimes
\alpha$ and
the proof is complete in this case.

\underbar{Case 2} Rank $\theta \le 1$ everywhere on $X$.  There exists a
collection of points $(p_1, \cdots, p_k)$ such that Ker $\theta$ extends
to a rank one subbundle of $T X \otimes \xi$, say $\alpha \otimes \eta$.  Since
$H^2 (X - \{ p_1, \cdots, p_k \}) \approx H^2 (X), c_1 ( \alpha)$ is
well-defined in $H^2 (X, \bbz)$.
Moreover, by the removing of singularities in codimension two we have $H^i (X,
{\cal O} )  \approx  H^i (X - \{ p_1 \cdots p_k \}, {\cal O})$, so from the
exact sequence $0 \to \bbz \to \cal O \to \cal O^\ast \to 1$
and the five-lemma we deduce that $H^1 (X, {\cal O}^\ast ) \approx  H^1 (X - \{
p_1, \cdots, p_k \},
\cal O^\ast)$, so $c_1 (\alpha)$ is in the image of $Pic (X)$ in $H^2 (X,
\bbz)$.
Let
$C$ be an irreducible curve of sufficiently high degree, which does not
meet $p_1, \cdots, p_k$.  Since $T X \otimes \ell \oplus \ell$ remains
$\theta$-stable
on $C$ [S1] we get $(c_1 (\alpha \otimes \ell) |_C, [C]) < 0$.  In view of
$H^{1,1} (X ) \cap H^2 (X_1 \bbq) \approx \bbq$ we can rewrite this as $(c_1
(\alpha), [ \ell] ) < - [ \ell]^2$.
Since $[\ell]$ generates $Pic (X) / tors$, this actually means
$(c_1 (\alpha), [ \ell]) \le
- 2 [ \ell]^2$.  Now, $c_1 (T X) = - 3 [ \ell]$, so $(c_1 ( T X / \alpha), [
\ell] ) \ge - [ \ell]^2$.  On $C$
we have an isomorphism
$$ \theta|_C : T X | \alpha \otimes \xi \to Im \theta \subset \eta |_C. $$
Hence $(c_1 (\text{Im} \theta), [C]) = (c_1 (T X / \alpha) + c_1 (\xi) ,
[C] ) \ge
(c_1 (\xi) - [ \ell], [ C])$.  Since, again, $[\ell]$ generates
$Pic (X) / tors$ and
$(c_1 (\xi) , [ \ell] > 0$ we get $(c_1 (\xi) - [ \ell], [ C ] ) \ge
0$.  This
contradicts the $\theta$-stability of $\xi \oplus \eta |_C$, because
$\text{Im} \theta |_C$ is $\theta$-invariant.  The proof is complete.

\centerline{References}

\item{[C]} K. Corlette,{\it Archimedian superrigidity and hyperbolic geometry},
Annals of Math., {\bf 135}, (1991), 165--182.

\item{[DM]} P. Deligne, G.Mostow,{\it Commensurabilities Among Lattices in
$PSU(1,n)$}, Annals of Mathematical Studies, {\bf 132}, Princeton Univ. Press,
1993.

\item{[GS]} M. Gromov, R. Schoen, {\it Harmonic maps into singular spaces and
$p$-adic superrigidity for latices in groups of rank one}, Publ.Math.IHES,{\bf
76} (1992),165--246.

\item{[JM]} D. Johnson and J. Millson, {\it Deformation spaces associated to
compact
hyperbolic manifolds}, in:Discrete Groups in Geometry and Analysis (papers in
honor of G.D.Mostow), Progress in Math.{\bf 67}, (1987), 48--106.

\item{[Re]} A. Reznikov,{\it  All regulators of flat bundles are torsion},
Annals of Math., (1995), to appear.

\item{[Rog]} J. Rogawski,{\it  Automorphic Representation of Unitary Groups in
Three Variables},
Sect. 15.3, Princeton Univ. Press, 1990.

\item{[S1]} C. Simpson, {\it Higgs bundles and local systems}, Publ. Math.
IHES 75 (1992), 5--95.

\item{[S2]} C. Simpson,{\it Moduli of representations of the fundamental group
of a smooth projective variety I}, Publ.Math.IHES {\bf 79} (1994),47--129.

\item{[W]} A. Weil, {\it Discrete subgroups of Lie groups I, II}, Annals of
Math.,
{\bf 72} (1960), 369--384, {\bf 75} (1962), 578--602.

\item{} Institute of Mathematics

\item{} Hebrew University

\item{} Givat Ram 91904

\item{} ISRAEL

\item{} email: simplex\@sunset.huji.ac.il
\bye